\documentclass{article}
\usepackage[all, knot]{xy}
\usepackage{amsbsy,amssymb,amscd,amsfonts,latexsym,amstext,delarray,
amsmath, diagrams}  \headsep=15pt
\def\hpi{\hat{\pi}}

\def\cX{{\cal X} }

\def\bX{ \bar{X}}

\def\dim{{ \mbox{dim} }}
\def\Spec{{ \mbox{Spec} }}

\def\ra{{ \rightarrow }}

\def\a{{ \alpha }}

\def\F{ {\mathbb F} }

\def\hra{{ \hookrightarrow }}

\def\C{{ \mathbb{C} }}

\def\bs{ \backslash}

\def\G{{ \Gamma }}
\def\Gal{{ \mbox{Gal} }}

\def\bQ{\bar{\Q}}

\def\bb{ \bar{b} }

\def\Z{{ \mathbb{Z}}}

\def\bq{\begin{quote}}
\def\eq{\end{quote}}

\newtheorem{obs}{Observation}
\def\Q{\mathbb{Q}}

\def\P{ {\bf P}}

\def\be{\begin{equation}}
\def\ee{\end{equation}}
\def\k{ \kappa}

\def\D{ \Delta}

\def\uet{U^{et}}

\def\hk{\hat{\k}}
\def\bx{ \bar{x}}
\def\loc{\mbox{loc}}
\def\bb{\bar{b}}

\title{Remark on fundamental groups and effective Diophantine methods for hyperbolic curves }
\author{Minhyong Kim}
\begin{document}
\maketitle
{\em Dedicated to the memory of Serge Lang}
\medskip

In a few earlier papers (\cite{kim1}, \cite{kim2}, \cite{KT}) attention was called
to the striking parallel between the ideas surrounding the well-known conjecture
of Birch and Swinnerton-Dyer for elliptic curves, and the mysterious
{\em section conjecture} of Grothendieck \cite{grothendieck} that concerns hyperbolic curves.
We wish to explain here some preliminary ideas
for `effective non-abelian descent' on hyperbolic curves equipped with
at least one rational point.
We again follow in an obvious manner the method of descent
on elliptic curves and, therefore, rely on conjectures.
In fact, the main point is to substitute the section conjecture for
the finiteness of the Shafarevich group. That is to say,
the input of the section conjecture is of the form
\bq
section conjecture $\Rightarrow$ termination of descent
\eq
At a number of different lectures on the topic of fundamental groups
and Diophantine geometry, the question was raised
about  the role of {\em surjectivity}
in the section conjecture as far as Diophantine
applications are concerned. This implication is intended as
something of a reply.

To {\em start} the descent, on the other hand, requires
the use of $p$-adic Hodge theory and the unipotent Albanese map.
In this process, in general, one another conjecture is unfortunately
needed. It could be, for example, the Bloch-Kato conjecture
on surjectivity of the $p$-adic Chern class map that has been
referred to in \cite{kim2}. In other words,
via the construction of Selmer varieties and
Albanese maps, one deduces an implication
\bq
Bloch-Kato conjecture $\Rightarrow$ beginning of descent
\eq
The main caveat here arises from lack of actual knowledge of
computational issues on the part of the author. To avoid misleading anyone about what is being achieved here,
we have in the following section separated out the questionable portions as hypotheses
[H] and [H'].
That is to say, the objects that mediate this process, namely
Galois cohomology groups/varieties and maps between them, seem in principle
to be computable.  But even to the algorithmically illiterate
perspective, it is obvious that actual computation
would be daunting to the point of impossibility given the technology of the
present day. Nevertheless, it is perhaps not entirely devoid of value
to point out  one direction of investigation in effective
methods, in the hope that even incompetent strategies may  eventually be
improved through the focussing  of sharper skills obviously available in the
community. Hence, the present paper.

One  point of some theoretical interest concerns
the comparison with `effective Mordell conjectures' in the
usual sense where upper bounds for heights are proposed.
If we fix a point $b$ on the curve and measure heights with respect
to the corresponding divisor, the height of another point
measures the distance from $b$ at all places. So an upper bound
for the height corresponds to a lower bound for the distance from $b$
at all places.
On the other hand, what the $p$-adic Hodge theory provides
(in principle) is a lower bound for the $p$-adic distance between
 all pairs of points at one place. This lower bound is exactly what is
 required to start the descent.

 Finally, we make the obvious point that the use of conjectures
 is probably not a serious obstacle from the computational
 perspective (that is, in comparison to the problem of
 feasibility). This is in the same spirit as the standard algorithms
 for  computing Mordell-Weil groups of elliptic
 curves where the BSD conjecture is employed with just a few
 misgivings \cite{cremona}.

 \section{Brief review}
 Here we will be intentionally brief, referring the reader to \cite{delignefg} and
 \cite{kim2} for a more thorough discussion.

Let $X/\Q$ be a proper smooth hyperbolic curve of genus $g$ with a point $b\in X(\Q)$
and let $S$
be the set of primes of bad reduction for $X$. In the following, we shall be a bit sloppy and
mostly omit separate notation for an integral model of $X$.
Choose a prime $p\notin S$ and
let $\uet=\pi^{et,\Q_p}_1(X,b)$  be the
$\Q_p$-unipotent \'etale fundamental group of
$\bX:=X\times_{\Spec(\Q)}\bQ$ and $\uet_n=(\uet)^n\bs \uet$
its quotient by the $n$-th level of the descending central
series normalized so that $(\uet)^1=\uet$. Let $\G$ be the
Galois group of $\bQ$ over $\Q$.
We defined (\cite{kim2}, \cite{KT})
the Selmer varieties
$$H^1_f(\G, \uet_n)$$
classifying $\G$-equivariant torsors for $\uet_n$
that are unramified at all places not in $\{p\}\cup S$
and crystalline at $p$.
($H^1_f(\G_T, \uet_n)$ in the notation of \cite{kim1} and \cite{kim2}.)
Recall the fundamental diagram (\cite{kim2}, end of section 2)
$$\begin{diagram}X(\Q) & \rInto &X(\Q_p)  &  &\\
\dTo^{\k^{et, glob}_n} & & \dTo^{\k^{et, loc}_{n} }& \rdTo^{\k^{dr/cr}_{n}} &\\
H^1_f(\G, U_n^{et}) &  \rTo^{\mbox{loc}_p} & H^1_f(G_p, U_n^{et}) & \rTo^{D} &U^{dr}_n/F^0\\
\end{diagram}$$
Here,
$H^1_f(G_p, U_n^{et})$ classifies $G_p:=\Gal(\bQ_p/\Q_p)$-equivariant torsors
for $\uet_n$ that are crystalline, while
$U^{dr}_n/F^0$ classifies compatible pairs
$T^{dr}_n\simeq T^{cr}_n$ of torsors for the De Rham and crystalline fundamental groups
$U^{dr}_n$  and $U^{cr}_n$ equipped with Hodge filtrations and Frobenius endomorphisms
compatible with the torsor structures.
The maps associate to each point $x\in X(\Q)$ the class of the torsor of paths from $b$ to $x$ in
the appropriate category. So
$$\k^{et,glob}_n(x)=[\pi^{et,\Q_p}_1(\bX;b,x)_n]$$ with $\G$-action,
$$\k^{et,loc}_n(x)=[\pi^{et,\Q_p}_1(\bX;b,x)_n]$$ with $G_p$-action,
and
$$\k^{dr/cr}_n(x)=[\pi^{dr}_1(X\otimes \Q_p; b,x)\simeq \pi^{cr}_1(Y; \bb,\bx)]$$
where $Y$ is the reduction mod $p$ of a smooth $\Z[1/S]$ model for $X$.

In contrast to this mass of notation, the section conjecture considers
just one map
$$\hk: X(\Q) \ra H^1(\G, \hpi_1(\bX,b))$$
that sends a point $x\in X(\Q)$
to the class of the pro-finite torsor of paths
$$\hpi_1(\bX;b,x)$$
with $\G$-action. It proposes that this map should be a bijection.
The injectivity is already known as a consequence of the Mordell-Weil
theorem for the Jacobian $J$ of $X$, while the surjectivity seems
to be a very deep problem. The question mentioned in the introduction arises
exactly because the injectivity appears, at first glance, to be more
relevant for finiteness than the surjectivity. The idea for
using the {\em bijectivity} seems to have been to create a tension between
the compact pro-finite topology of $
H^1(\G, \hpi_1(\bX,b))$ and the `discrete nature' of $X(\Q)$.
At present it is unclear how this intuition is to be realized.
But, as mentioned, when the finiteness is obtained through a different
approach, we wish to explain the use of the surjectivity for
{\em finding} the full set of points.

Using the exact sequence
$$0\ra U^{n+1}\bs U^n \ra U_{n+1} \ra U_n \ra 0$$
for each of the fundamental groups, the
global Selmer variety is fibered according to the sequence
$$0\ra H^1_f(\G,(\uet)^{n+1}\bs (\uet)^n) \ra H^1_f(\G, \uet_{n+1}) \ra H^1_f(\G, \uet_n)$$
which means that the kernel acts on the variety in the
middle with orbit space a subset of the third object.
If we denote by $r_n$ the dimension of $U_n$, there is a recursive formula \cite{labute}
$$\Sigma_{i|n}ir_i=(g+\sqrt{g^2-1})^n+(g-\sqrt{g^2-1})^n$$
which implies in particular that
$$r_n\approx (g+\sqrt{g^2-1})^n/n.$$

The global Selmer variety has its dimension controlled by the
Euler characteristic formula for the cohomology of the group
$G_T=\Gal (\Q_T/\Q)$, where $T=S\cup \{p\}$ and $\Q_T$ is the maximal
extension of $\Q$ unramified outside $T$ (\cite{kim2}, section 3). It reads
$$H^1(\G_T, (\uet)^{n+1}\bs (\uet)^n) - H^2(\G_T, (\uet)^{n+1})\bs (\uet)^n)
= [(\uet)^{n+1}\bs (\uet)^n]^{-}$$
where the minus in the superscript refers to the sign for the action of
complex conjugation. The dimension of this
minus part can be estimated as follows.
The action of complex conjugation on the
\'etale fundamental group is compatible
with the  its action on  the Betti realization
of the motivic fundamental group \cite{delignefg} according to which
$$
(U^B)^{n+1}\bs (U^B)^n$$
has a pure Hodge structure of weight $n$.
So when $n$ is odd, we get
$$\dim  [(\uet)^{n+1}\bs (\uet)^n]^{-}=r_n/2$$
But when
$n=2m$ is even, there is the contribution from
the $(m,m)$ component to consider, which can be complicated.
This $(m,m)$ component is a quotient of the $(m,m)$-part
of $$H_1(X(\C), \C)^{\otimes 2m}$$ which has dimension
$\binom{2m}{m} g^{2m}$. So for simplicity, we will just use the tautological
estimate
$$\dim  [(\uet)^{n+1}\bs (\uet)^n]^{-}\leq r_n$$
for $n$ even.

In \cite{kim2}, section 3, we analyzed already the use of the corresponding Shafarevich groups
$$Sha^2((\uet)^{n+1}\bs (\uet)^n)):=
Ker [H^2(\G_T, (\uet)^{n+1}\bs (\uet)^n)\ra
\oplus_{v\in T} H^2(\G_v, (\uet)^{n+1}\bs (\uet)^n)]$$
which is dual to
$$Sha^1(((\uet)^{n+1}\bs (\uet)^n)^*(1)):=
Ker [H^1(\G_T, ((\uet)^{n+1}\bs (\uet)^n)^*(1))\ra
\oplus_{v\in T} H^1(\G_v, ((\uet)^{n+1}\bs (\uet)^n)^*(1))]$$
There is a Chern class map \cite{BK}
$$ch_{n,1}: K^{(1)}_{2-n-1}(X^n)\otimes \Q_p\ra H^1(\Gal(\bar{\Q}/\Q), H^n(\bX^n, \Q_p(1)))$$
for $n\neq 1$
whose image lies in a `geometric' subspace
$$H^1_g(\Gal(\bar{\Q}/\Q), H^n(\bX^n, \Q_p(1)))$$
that contains
$$Sha^1(H^n(\bX^n, \Q_p(1)))$$
In fact,
$$Sha^1([(\uet)^{n+1}\bs (\uet)^n)]^*(1))$$
is a subspace of $Sha^1(H^n(\bX^n, \Q_p(1)))$ because
the representation
$(\uet)^{n+1}\bs (\uet)^n$ is a direct summand of
$H^{et}_1(\bX, \Q_p)^{\otimes n}$
which, in turn, is a direct summand of
$(H^n(\bX^n, \Q_p))^*$.
But Bloch and Kato conjecture that
$$ch_{n,1}: K^{(1)}_{2-n-1}(X^n)\otimes \Q_p \ra H^1_g(\Gal(\bar{\Q}/\Q), H^n(\bX^n, \Q_p(1)))$$
is an isomorphism.
Thus, when $n\geq 2$, we get
$$ Sha^1([(\uet)^{n+1}\bs (\uet)^n)]^*(1))=0$$

We recall the explicit bound for the local $H^2$ (\cite{kim2}, section 3).
For $v\neq p$, we have
$$\dim H^2(G_v, (\uet)^{n+1}\bs (\uet)^n) \leq n
g^{n}+ \frac{n(n-1)}{2} (2g-2)^2 g^{n-2}$$
while
$$\dim H^2(G_p, (\uet)^{n+1}\bs (\uet)^n) \leq n g^{n}$$

Finally, as regards the contribution of the
Hodge filtration, we saw in loc. cit. that
$$F^0((U^{dr})^{n+1}\bs (U^{dr})^{n+1})\leq g^n$$
so that
$$\dim(U^{dr})^{n+1}\bs (U^{dr})^{n+1}/F^0\geq r_n-g^n\approx (g+\sqrt{g^2-1})^n/n-g^n$$

\section{Beginning the descent}
Since it costs very little extra work to define,
we will in fact consider the refined Selmer variety
$$H^1_{f,0}(\G, \uet_n)\subset H^1_{f}(\G, \uet_n)$$
consisting of classes whose images in
$$H^1_{f}(\G, \uet_2)$$ go
to zero under all localization maps
$$H^1_{f}(\G, \uet_2) \stackrel{\loc_v}{\ra} H^1(G_v, \uet_2)$$
for $v\neq p$.
As explained in \cite{KT},
the image of $X(\Q)$ under $\k^{et,glob}_n$
lies in
$H^1_{f,0}(\G, \uet_n)$.
From the estimates of the previous
section, it is obvious
that
\smallskip

\begin{em}
assuming the Bloch-Kato conjecture, there is an effectively computable $t$ such that
$$\dim H^1_{f,0}(\G, \uet_n)< \dim U^{dr}_n/F^0$$
for $n\geq t$.
\end{em}
\smallskip

Of course the computation starts out with an estimate for
$\dim H^1_{f,0}(\G, \uet_2)$
which according to the usual BSD is the same as the
Mordell-Weil rank of $J$.
After that the dimension of
$\dim H^1_{f,0}(\G, \uet_n)$ grows as a function of $n$ with
an explicit upper bound  while the dimension of $U^{dr}_n/F^0$
grows with an explicit (and eventually bigger) lower bound.
Written out, the estimate for growth looks like
$$ \dim H^1_{f,0}(\G, \uet_{2n+1})
\leq \dim H^1_{f,0}(\G, \uet_{2n})+r_{2n}/2+|S|[(2n)
g^{2n}+ \frac{(2n)(2n-1)}{2} (2g-2)^2 g^{2n-2}] +(2n)g^{2n}$$
and
$$\dim H^1_{f,0}(\G, \uet_{2n+2})
\leq \dim H^1_{f,0}(\G, \uet_{2n+1})+r_{2n+1}+|S|[(2n+1)
g^{2n+1}+ \frac{(2n+1)(2n)}{2} (2g-2)^2 g^{2n-1}] +(2n+1)g^{2n+1}$$
while
$$\dim U_{n+1}\geq \dim U_n+r_n-g^n$$
We eventually get an inequality in the right direction because of the asymptotic behavior of $r_n$.
In this regard, note that $g+\sqrt{g^2-1}>g$ for $g\geq 2$.

As a consequence of the discrepancy in dimension,
the image of
$$D\circ \loc_p:  H^1_{f,0}(\G, \uet_t)\ra U^{dr}_t/F^0$$
is {\em not} Zariski dense.
In contrast to difficult sets like $X(\Q)$, the classifying spaces for torsors and the
maps between them are algebro-geometric objects
which can be computed in principle. This should work in the manner of computations with the usual method of Chabauty
as appears, for example, in  \cite{flynn} (cf. the discussion of $\theta$ in the introduction).
In case this is not convincing, we will adopt it as an additional
hypothesis:

\bq
$[H]$: The map
$$D\circ \loc_p: H^1_{f,0}(\G, \uet_t)\ra  U^{dr}_t/F^0$$
can be computed.
\eq
The end result of this is
that assuming B-K and [H], we can find  an algebraic
function $\a$ on $U^{dr}_t/F^0$, that vanishes on
the image of $H^1_{f,0}(\G, \uet_t)$.
Now, when we restrict $\a$ to
$X(\Q_p)$ it becomes a linear combination of $p$-adic iterated
integrals. To elaborate on this point a little more, recall (\cite{kim2}, section 1) the description
of the coordinate ring of the De Rham fundamental group $U^{dr,0}$
for an affine curve $X^0$ obtained by deleting some rational
divisor from $X$. In this case, when we choose a collection
$a_1, a_2, \ldots, a_k$ of  algebraic differential forms
on $X^0$ inducing a basis of $H^1_{dr}(X^0)$,
the coordinate ring of
$U^{dr,0}$ has the form
$$\Q_p<a_w>,$$
the $\Q_p$ vector space generated by symbols
$a_w$, one for each finite sequence $w$ of numbers from $\{1, 2, \ldots, k\}$.
Furthermore, on
$X^0(\Z_p)$, there is a lifting (depending on the previous choice of
basis)

$$\begin{diagram}
 & & U^{dr,0}_t \\
 & \ruTo & \dTo \\
 X^0(\Z_p) & \rTo & U^{dr,0}_t/F^0
 \end{diagram}$$
such that the restriction of $a_w$ for
$w=(i_1, i_2, \ldots, i_l)$ to $X^0(\Z_p)$ has the form
$$a_w(z)= \int_b^{z} a_{i_1}a_{i_2}\cdots a_{i_l}$$
Also, there is a functorial map
$$U^{dr,0}_t \ra U^{dr}_t$$
compatible with the Hodge filtration so that the function
$\a$ on $U^{dr}_t/F^0$ can be lifted to $U^{dr,0}_t$. That is to say, one can construct a diagram
$$\begin{diagram}
 & & U^{dr,0}_t \\
 & \ruTo & \dTo \\
 X^0(\Z_p) & \rTo & U^{dr}_t/F^0
 \end{diagram}$$
 enabling us to compute the restriction of $\a$ to
 $X^0(\Z_p)$ in terms of the $a_w$.
  The idea would be to carry  this process out for
two separate affine $X^0$ so as to cover $X(\Z_p)$ and
then to express $\a$  in terms of iterated integrals on each affine open set.
 Of course, the problem of explicitly computing
the local liftings is also a daunting task,
although possible in theory. The author makes no pretense of
knowing, as yet, how to reduce this to a tractable
process. Perhaps it is safer to state it also explicitly as a hypothesis:
\bq
[H']: The map
$$\begin{diagram} U^{dr,0}_t & \rTo & U^{dr}_t/F^0 \\
\end{diagram} $$
can be computed.
\eq
Choose a representative $y\in X(\Q_p)$ for
each point in $Y(\F_p)$ ($=X$ mod $p$) and a coordinate $z_y$ centered at $y$. We must then approximate the zeros of
$\a$ on $X(\Q_p)$ by expressing it as a power series in the
$z_y$. This needs to be carried out to a sufficiently high degree of accuracy
so that we can find an $M$ and a finite collection $y_i\in X(\Q_p)$  for which
$$]y_i[_M:=\{x\in X(\Q_p)| z_{y_i}(x) \leq p^{-M}\}$$
contains at most one zero of $\a$. That is to say, we need to separate the
zeros of $\a$ modulo $p^M$.
Note that even at this point, since all expressions
will be approximate, there would be no way to determine which of the $y_i$
relate to actual points of $X(\Q)$, even  though an upper bound for
the {\em number} of points may be available, as was emphasized by Coleman \cite{coleman}.
In fact,  the process of separating the points using small disks already seems to
occur, at least implicitly, in the method of Coleman-Chabauty. In the next section
we will see how to combine that separation with the section conjecture.

We summarize the preceding passages as follows:

\begin{obs} Assuming the Bloch-Kato conjecture and  the hypotheses [H]
and [H'], there is an effectively
computable $M$ such that the
map
$$X(\Q) \hra X(\Q_p) \ra X(\Z/p^M)$$
is injective.
\end{obs}

In our view, this statement is one rather essential justification for studying the
Selmer varieties and unipotent Albanese maps. That is, Faltings' theorem
as it stands does not seem to give, even in principle, a way of getting
at this sort of effectivity. To belabor the obvious, the point
is that the map
$$X(\Q) \ra X(\Q_p)$$ is not a priori (i.e., before finding $X(\Q)$)
computable  even in
principle, while
$$H^1_{f,0}(\G, \uet_t)\ra  U^{dr}_t/F^0$$
is.

When we embed $X(\Q)$ inside
$J(\Q)$ using the base-point $b$, we see then
that we have an injection
$$X(\Q) \hra J(\Z/p^M)$$
But the kernel of the reduction map
$$J(\Q) \ra J(\Z/p^M)$$
is of finite index, and hence, contains
$ NJ(\Q)$ for some $N$.
 So finally, we arrive at an effectively computable $N$ such that
 $$X(\Q)\ra J(\Q) \ra  J(\Q)/NJ(\Q)$$
 is injective.
 Let $T_0$ be $S$ together with the set of primes dividing $N$ and $\G_{T_0}$
 the fundamental group of $\Spec (\Z[1/T_0])$ with base-point
 given by
 $\Z[1/T_0]\hra \Q \hra \bQ$. Then
 we get an injection
 $$X(\Q) \hra H^1(\G_{T_0}, J[N])$$
 allowing us to begin descent.

\section{Non-abelian descent and its termination}
Once we have the final conclusion of the previous section, we can dispense
entirely with the unipotent machinery and start to deal with the pro-finite formalism.
There are many ways to construct a co-final system for
$$\D:=\hpi_1(\bX,b)$$
of which we will use one described in a letter from Deligne to
Thakur \cite{delignethakur}.
Let $K_n\subset \D$ be the intersection of all open subgroups of
index $\leq n$. It is a characteristic subgroup, and hence,
we can form the quotient $\D(n):=\D/K_n$.
The order of this quotient has all prime divisors $\leq n$.
Let $\G_n$ denote the fundamental group of
$\Spec(\Z[1/n!])$. We also denote by
$\pi(n)$ the quotient of $\hpi_1(X,b)$
by $K_n$, a group that fits into  the exact sequence
$$0\ra \D(n) \ra \pi(n) \ra \G \ra 0.$$
For $n$ larger than any prime in $S$, there is a pull-back diagram (\cite{wildeshaus}, proof of theorem 2.8)
$$\begin{diagram}
0 & \rTo & \D(n) & \rTo &\pi(n) & \rTo & \G & \rTo& 0 \\
 & & \dTo^{=} & & \dTo & & \dTo & & \\
0 & \rTo & \D(n) & \rTo & \hpi_1(\cX_n)/K_n& \rTo & \G_n & \rTo& 0
\end{diagram}$$
where
$\cX_n$ is a proper smooth  model for $X$ over $\Spec(\Z[1/n!])$.
Therefore, we see that
any point $x\in X(\Q)$ defines a class in
$$H^1(\G_n, \D(n))$$ and that we have a commutative diagram
\be\begin{diagram}
X(\Q) & \rInto^{\hk} &H^1(\G, \D) \\
\dTo & & \dTo \\
H^1(\G_n, \D(n)) & \rInto & H^1(\G, \D(n))
\end{diagram}\ee
There is  a  sequence of subsets containing $X(\Q)$,
$$
H^1(\G, \D)_i \subset H^1(\G, \D)$$
consisting of those classes
whose projection to
$H^1(\G, \D(i))$ lie in the image of
$$H^1(\G_i, \D(i) )\hra H^1(\G, \D(i))$$
Let $n_0$ be larger than the primes in $T_0$. Then we have
 diagrams
\be\begin{diagram}
 & & H^1(\G, \D)_i & \rInto & H^1(\G, \D)\\
 & &       \dTo    &       & \dTo\\
 & & H^1(\G_i, \D(i)) & \rInto & H^1(\G, \D(i))\\
 & & \dTo & & \\
 H^1(\G_{T_0}, J[N])& \rInto & H^1(\G_i, J[N])& &
 \end{diagram}\ee
 for $i\geq n_0$.
 Using this, we can define a decreasing sequence
 of subsets
 $$ H^1(\G_{T_0}, J[N])_n \subset  H^1(\G_{T_0}, J[N])$$
for $n\geq n_0$ consisting of those classes whose images in
$H^1(\G_i, J[N])$ lift to
$H^1(\G_i, \D(i))$ for all $n_0 \leq i\leq n$.
For $n\geq n_0$, we also have a commutative diagram
\be \begin{diagram}
X(\Q) & \rInto & H^1(\G_n, \D(n)) \\
\dInto & & \dTo \\
H^1(\G_{T_0}, J[N]) & \rInto & H^1(\G_n, J[N])
\end{diagram}\ee

Meanwhile, there is  an increasing sequence
$$X(\Q)_n\subset X(\Q) \subset H^1(\G_{T_0},J[N])$$
consisting of the points with height (in some projective embedding)
$\leq n$.
We visualize the situation using the sort of diagram familiar from the
arithmetic theory of elliptic curves:
$$
\cdots X(\Q)_n \subset X(\Q)_{n+1} \subset \cdots \subset H^1(\G_{T_0}, J[N])_{m+1} \subset
H^1(\G_{T_0}, J[N])_m \subset \cdots \subset H^1(\G_{T_0}, J[N])$$

\begin{obs}
The section conjecture implies that
$$X(\Q)_n=H^1(\G_{T_0}, J[N])_m $$
for $n, m$ sufficiently large. At this point,
$X(\Q)=X(\Q)_n$.
\end{obs}

That is to say,  we know when to stop searching.
The simple proof is written out just to make sure the author is not
confused.

{\em Proof.}
Assume the section conjecture. Then by diagrams (1) and (2),
we have
$$H^1(\G,\D)_i=H^1(\G,\D)$$
for all $i$ and we actually have maps
$$H^1(\G,\D) \ra H^1(\G_i, \D(i))$$
for each $i$. Furthermore since $H^1(\G, \D)$ is finite
(which follows either from Faltings theorem or the
reproof assuming Bloch-Kato from the previous section)
we have
$$H^1(\G,\D) \simeq H^1(\G_i, \D(i))$$
for $i$ sufficiently large. So if $c\in H^1(\G_{T_0}, J[N])$
is not in $X(\Q)$, then $c\notin H^1(\G_{T_0}, J[N])_m$
for some $m$. Thus, eventually,
$X(\Q)=H^1(\G_{T_0}, J[N])_m$. Of course eventually
$X(\Q)_n=X(\Q)$.
Now suppose
$$X(\Q)_n=H^1(\G_{T_0}, J[N])_m$$ at any point.
Then classes not in $H^1(\G_{T_0}, J[N])_m$
cannot lift to $H^1(\G,\D)_m$. And hence,
they are not in $X(\Q)$. That is to say,
$X(\Q)_n=X(\Q)$.
$\Box$

All the cohomology sets occurring in the  argument
are finite and thereby have the nature of being computable
through explicit Galois theory.
As mentioned in the introduction, the actual implementation of
such an algorithm is obviously an entirely different matter.
\begin{flushleft}
{\bf Acknowledgements:}

-The author was supported in part by
a grant from the National Science Foundation and a visiting professorship
at RIMS.

-He is grateful to Kazuya Kato, Shinichi Mochizuki, and Akio
Tamagawa for a continuing stream of discussions on topics related
to this paper, and for their generous hospitality during the Fall
of 2006.
\smallskip

\end{flushleft}

{\footnotesize DEPARTMENT OF MATHEMATICS, PURDUE UNIVERSITY,
  WEST LAFAYETTE, INDIANA 47907  and DEPARTMENT OF MATHEMATICS, UNIVERSITY OF ARIZONA,
TUCSON, AZ 85721, U.S.A. }

{\footnotesize EMAIL: kimm@math.purdue.edu}

 \end{document}